\documentclass[reqno, a4paper, 12pt]{amsart} %
\usepackage{epsfig}   
\usepackage{ifpdf} 
\usepackage{amssymb,amsmath}         
\usepackage{graphicx}   
\usepackage{float}      

\usepackage[paperwidth=8.3in,paperheight=11.7in,
  left=0.7in,right=0.7in,top=1.0in,bottom=1.0in]{geometry}

\begin{document}

\makeatletter
\def\rank{\mathop{\operator@font rank}\nolimits}
\def\det{\mathop{\operator@font det}\nolimits}
\makeatother

\newtheorem{thm}{Theorem}[section]
\newtheorem{ex}[thm]{Example}
\newtheorem{lem}[thm]{Lemma}
\newtheorem{rmk}[thm]{Remark}
\newtheorem{defi}[thm]{Definition}
\newtheorem{cor}[thm]{Corollary}


\title[asymptotic stabilization]{Design of the state feedback-based feed-forward controller asymptotically stabilizing the overhead crane at the desired end position}

\author{Robert Vrabel} 
\address{Robert Vrabel, Slovak University of Technology in Bratislava, Faculty of Materials Science and Technology, Institute of Applied Informatics, Automation and Mechatronics,  Bottova~25,  917 01 Trnava, Slovakia}
\email{robert.vrabel@stuba.sk}

\date{{\bf\today}}

\begin{abstract}
The problem of feed-forward control of overhead crane system is discussed. By combining the Kalman's controllability theory and Hartman-Grobman theorem from dynamical system theory, a linear, continuous state feedback-based feed-forward controller that stabilizes the crane system at the desired end position of payload is designed. The efficacy of proposed controller is demonstrated by comparing the simulation experiment results for overhead crane with/without time-varying length of hoisting rope.
\end{abstract}


\maketitle

\section[Introduction]{Introduction}

For overhead crane control, it is required that the trolley should reach the desired location as fast as possible while the payload swing should be kept as little as possible during the transferring process. However, it is extremely challenging to achieve these goals simultaneously owing to the underactuated characteristics of the crane system, more specifically, underactuated with respect to the load sway dynamics, which makes the linear part of overhead crane system uncontrollable in the sense of Kalman's control theory by a continuous control law.  Due to this reason, the development of efficient control schemes for overhead cranes has attracted wide attention from the control community. 
For example, in \cite{Park} a nonlinear control law for container cranes with load hoisting using the feedback linearization technique and the decoupling strategy of swing angle-dynamics from trolley movement- and varying rope length- dynamics by Lyapunov function approach was investigated. 
The authors in \cite{Tagawa_et_al} and  \cite{Tagawa2_et_al} developed a sensorless vibration control system for overhead crane system with varying wire length by using simulation-based control technique. 

In the paper \cite{Yang} an adaptive nonlinear coupling control law has been presented for the motion control of overhead crane with constant rope length and without considering the mass moment of inertia of the load. By utilizing a Lyapunov-based stability analysis, the authors achieved asymptotic tracking of the crane position and stabilization of payload sway angle of an overhead crane. 

A sliding mode controller composed of approximated control and switching action was designed in \cite{Tuan} for simultaneously combining control of cargo lifting, trolley moving, and cargo swing vanishing. 

In \cite{Tomczyk}, the problems of load operation and positioning under different wind disturbances by using dynamic model with a state simulator is discussed.

In all of these mentioned papers, but also in others, see, e.g. \cite{Sorensen}, \cite[p.~273]{Zhang} and the reference therein, the control of sway-angle of payload is not considered and, therefore, this angle is assumed to be small what allows some approximations and truncations in the reference model (Remark~\ref{approx}). Also recall that for the mathematical models without considering a mass moment of inertia of the payload, what is usual for the simplified models of crane systems, the model is singular for rope length near zero what makes the system practically unanalyzable from the point of view of engineering practice.  

So, the purpose of this article is to show that by adding the control force for load sway damping (in the case of varying length of rope), the linear part of system becomes state controllable and the original overhead crane system  asymptotically stabilizable at the required end position by linear and continuous state feedback and, moreover, the exact formulas for the control forces are given. The efficacy of proposed controller is demonstrated by comparing the computer simulation experiment results for overhead crane with and without varying length of hoisting rope, respectively.  

\begin{figure}
\ifpdf\vglue-4cm\includegraphics[trim={5cm 18cm 5cm 0},clip]{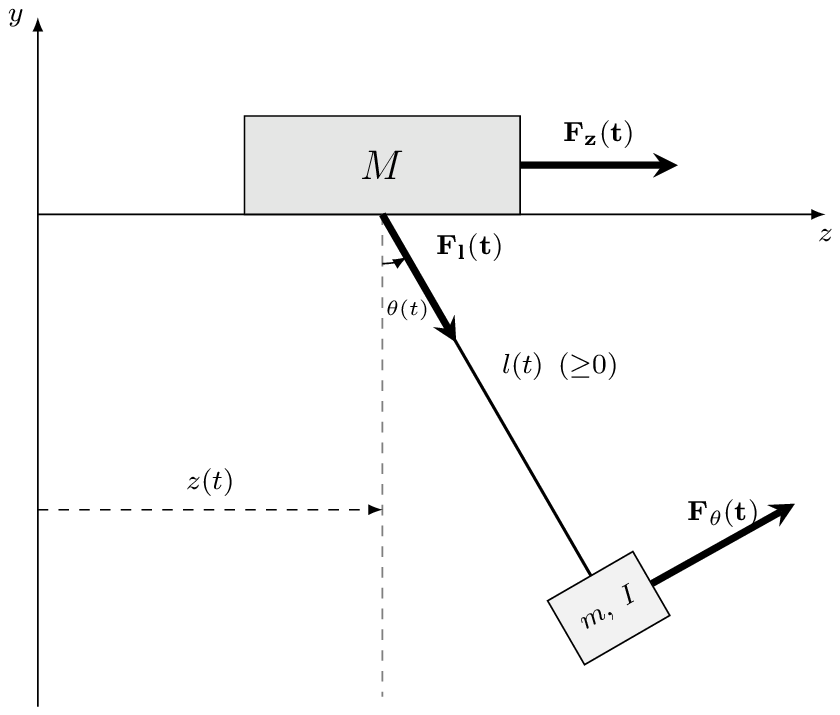}\else\includegraphics{crane_system.eps}\fi
\caption{Schematic diagram of the 2-D overhead crane system. The coordinates of payload are $z_p(t)=z(t)+l(t)\sin\theta(t),$ $y_p(t)=-l(t)\cos\theta(t)$} 
\label{fig:M1} 
\end{figure}

Anti-swing control of automatic overhead crane system required to transfer the payload without causing excessive swing at the end position. With a fully-automated overhead crane, the operators make the settings, and the crane automatically takes care of repetitive or difficult actions. This is especially useful in demanding and hazardous environments. Also, automated overhead cranes can reduce labor costs, track inventory, optimize storage, reduce damage, increase productivity and reduce the capital expense. Most of the proposed anti-swing controls use feedbacks that require two sensors to measure the trolley position and swing angle. However, installing swing angle sensor on a real overhead crane is often troublesome and also more costly. Moreover, vibration measurement sensors are required and this causes faults of the crane systems especially in severe industrial environments (\cite{Tagawa_et_al}).  Our approach is based on the feed-forward design of the control forces $F_z,$ $F_l$ and $F_{\theta},$ see Fig.~\ref{fig:M1}, applied to the crane system using simulation-based control strategy and that will stabilize the system at the desired and in advance known end position of the payload.

\section{Crane system dynamics. Reference model}

Basically, an overhead crane is made up of a trolley (cart) moving along a horizontal axis with a load hung from a flexible rope. The Fig.~\ref{fig:M1} shows the swing motion of the load caused by trolley movement, in which $z$ is the trolley moving direction, $y$ is the vertical direction,  $\theta(t)$ is
the sway angle of the load,  $z(t)$ is the position of the trolley, $l(t)$ is the hoist rope length. The symbols $F_z,$ $F_l$ and $F_{\theta}$ denote the control forces applied to the trolley in the  $z$-direction, to the payload in the  $l$-direction and in the direction perpendicular to $l$-direction, respectively.

The following assumptions are made throughout the paper: 
\begin{itemize}
\item[i)] The payload and trolley are connected by a massless, rigid rope, that is, a pendulum motion of the load is considered; 
\item[ii)] The trolley moves in the $z$-direction; 
\item[iii)] The payload moves on the $z-y$ surface;
\item[iv)] All frictional elements in the trolley and hoist motions can be neglected; 
\item[v)]  the rope elongation is negligible. 
\end{itemize}
The information on the sway angle, sway angular velocity, trolley displacement and velocity, hoisting rope length
and its time rate of change are not assumed to be known. So, this paper presents sensor less anti-swing control strategy for automatic overhead crane.

Using Lagrangian method, the Lagrangian equation with respective to the
generalized coordinate $q_i$ can be obtained as
\begin{equation}\label{Lagrange_eq}
\frac{d}{dt}\bigg(\frac{\partial L}{\partial\dot q_i}\bigg)-\frac{\partial L}{\partial q_i}=F_i,
\end{equation}
where $i = 1, 2, 3,$ $L = KE - PE$ ($KE$ means the system kinetic energy and $PE$ denotes the
system potential energy), $q_i$'s are the generalized coordinates, here $q_1,$ $g_2$ and $q_3$ indicate $z,$ $l$ and $\theta,$  respectively, and $F_i$ represents nonconservative generalized force ($F_z,$ $F_l$ and $F_{\theta}$) associated with those coordinates.
For the system under consideration the total kinetic energy and potential energy is 
\[
KE=\frac12(M+m){\dot z}^2+\frac12m{\dot l}^2+\frac12m(l\dot\theta)^2+m\dot z(\dot l\sin\theta+l\dot\theta\cos\theta)+\frac12I{\dot\theta}^2
\]
and
\[
PE=-glm\cos\theta,
\]
the potential energy of the trolley subsystem is kept unchanged. Here, $m$ is the payload mass, $M$ is the mass of the trolley with the hoist system,  $I$ is the mass moment of inertia of the
payload, and $g$ is the gravitational acceleration.

Using the Lagrangian equation (\ref{Lagrange_eq}), one obtain the following relations between the generalized coordinates $z,$ $l$ and $\theta:$
\begin{equation}\label{eq:z}
(M+m)\ddot z+m\ddot l\sin\theta+2m\dot l\dot\theta\cos\theta+lm\ddot\theta\cos\theta-lm{\dot\theta}^2\sin\theta=F_z,
\end{equation}
\begin{equation}\label{eq:l} 
m\ddot z\sin\theta+m\ddot l-lm{\dot\theta}^2-gm\cos\theta=F_l,
\end{equation}
and
\begin{equation}\label{eq:theta}
lm\ddot z\cos\theta+(ml^2+I)\ddot\theta+2lm\dot l\dot\theta+glm\sin\theta=F_{\theta}.
\end{equation}
Let us  introduce the state variables $x_1=:z,$ $x_2=:l,$ $x_3=:\theta,$ $x_4=:\dot z,$ $x_5=:\dot l,$ and $x_6=:\dot \theta.$ Now solving the equations (\ref{eq:z}), (\ref{eq:l}) and (\ref{eq:theta}) with regard the variables $\dot x_4,$ $\dot x_5$ and $\dot x_6,$ one obtain the linear system of equations
\[
\left(
\begin{array}{ccc} 
M+m & m\sin\left(x_{3}\right) & mx_{2}\cos\left(x_{3}\right)\\ 
m\sin\left(x_{3}\right) & m & 0\\ 
mx_{2}\cos\left(x_{3}\right) & 0 & m{x^2_{2}}+I 
\end{array}
\right)
\left(
\begin{array}{c}
\dot x_4\\
\dot x_5\\
\dot x_6
\end{array}
\right)
\]
\[
=\left(
\begin{array}{c} 
F_{z}+mx_{2}{x^2_{6}}\sin\left(x_{3}\right)-2mx_{5}x_{6}\cos\left(x_{3}\right)\\ 
F_{l}+mx_{2}{x^2_{6}}+gm\cos\left(x_{3}\right)\\
 F_{\theta}-2mx_{2}x_{5}x_{6}-gmx_{2}\sin\left(x_{3}\right) 
\end{array}
\right)
\]
and its solution
\[
\dot x_4= -\frac{1}{Mm{x^2_{2}}+Im{\cos^2\left(x_{3}\right)}+IM}
\]
\[
\times\Bigg[\cos\left(x_{3}\right)\left(F_{\theta}mx_{2}+2Imx_{5}x_{6}\right)
\]
\[
+\sin\left(x_{3}\right)\left(F_{l}m{x^2_{2}}+F_{l}I+Igm\cos\left(x_{3}\right)\right)
\]
\[
-F_{z}I-F_{z}m{x^2_{2}}\Bigg],
\]
\[
\dot x_5=\frac{1}{m\left(Mm{x^2_{2}}+Im{\cos^2\left(x_{3}\right)}+IM\right)}
\]
\[
\times\Bigg[F_{l}IM+F_{l}m^2{x^2_{2}}+F_{l}Im+Mm^2{x^3_{2}}{x^2_{6}}
\]
\[
+\frac{1}{2}F_{\theta}m^2x_{2}\sin\left(2x_{3}\right)-F_{z}m^2{x^2_{2}}\sin\left(x_{3}\right)+F_{l}Mm{x^2_{2}}
\]
\[
-F_{z}Im\sin\left(x_{3}\right)-F_{l}m^2{x^2_{2}}{\cos^2\left(x_{3}\right)}+Igm^2\cos\left(x_{3}\right)+Im^2x_{2}{x^2_{6}}{\cos^2\left(x_{3}\right)}
\]
\[
+Mgm^2{x^2_{2}}\cos\left(x_{3}\right)+IMgm\cos\left(x_{3}\right)+Im^2x_{5}x_{6}\sin\left(2x_{3}\right)+IMmx_{2}{x^2_{6}}\Bigg],
\]
and
\[
\dot x_6=\frac{1}{Mm{x^2_{2}}+Im{\cos^2\left(x_{3}\right)}+IM}
\]
\[
\times\Bigg[F_{\theta}M+F_{\theta}m{\cos^2\left(x_{3}\right)}-F_{z}mx_{2}\cos\left(x_{3}\right)+F_{l}mx_{2}\cos\left(x_{3}\right)\sin\left(x_{3}\right)
\]
\[
-2Mmx_{2}x_{5}x_{6}-Mgmx_{2}\sin\left(x_{3}\right)\Bigg],
\]
where $F_{z},F_l,F_\theta$ represent the control forces. 

Now let us substitute  $(M+m)u_1$ instead of $F_{z},$ $mu_2-gm\cos\left(x_{3}\right)$ instead of $F_{l}$ (here the term $-gm\cos\left(x_{3}\right)$ compensates the weight of load) and $Iu_3$ instead of $F_{\theta},$  where $u_1, u_2, u_3$ are new control variables. The reference mathematical model for design a feed-forward controller takes the form of control system  $\dot x =G(x,u),$ $G=(G_1,\dots,G_6)^T$ with
\[
G_1=x_4,
\]
\[
G_2=x_5,
\] 
\[
G_3=x_6,
\] 
\[
G_4= -\frac{1}{Mm{x^2_{2}}+Im{\cos^2\left(x_{3}\right)}+IM}
\]
\[
\times\Bigg[\cos\left(x_{3}\right)\left(Imu_{3}x_{2}+2Imx_{5}x_{6}\right)
\]
\[
+\sin\left(x_{3}\right)\left(m^2u_{2}{x^2_{2}}+Imu_{2}-gm^2{x^2_{2}}\cos\left(x_{3}\right)\right)
\]
\[
-Iu_{1}\left(M+m\right)-m\left(M+m\right)u_{1}{x^2_{2}}\Bigg],
\] 
\[
G_5=\frac{1}{{Mm{x^2_{2}}+Im{\cos^2\left(x_{3}\right)}+IM}}
\]
\[
\times\Bigg[IMu_{2}+m^2u_{2}{x^2_{2}}+Imu_{2}+IMx_{2}{x^2_{6}}
\]
\[
-gm^2{x^2_{2}}\cos\left(x_{3}\right)-IMu_{1}\sin\left(x_{3}\right)-m^2u_{1}{x^2_{2}}\sin\left(x_{3}\right)+Mmu_{2}{x^2_{2}}
\]
\[
-Imu_{1}\sin\left(x_{3}\right)+gm^2{x^2_{2}}{\cos^3\left(x_{3}\right)}-m^2u_{2}{x^2_{2}}{\cos^2\left(x_{3}\right)}
\]
\[
+Mm{x^3_{2}}{x^2_{6}}+\frac12 Imu_{3}x_{2}\sin\left(2x_{3}\right)+Imx_{5}x_{6}\sin\left(2x_{3}\right)
\]
\[
-Mmu_{1}{x^2_{2}}\sin\left(x_{3}\right)+Imx_{2}{x^2_{6}}{\cos^2\left(x_{3}\right)}\Bigg],
\]
and
\[
G_6= \frac{1}{Mm{x^2_{2}}+Im{\cos^2\left(x_{3}\right)}+IM}
\]
\[
\times\Bigg[IMu_{3}-x_{2}\left[2Mmx_{5}x_{6}+m\left(M+m\right)u_{1}\cos\left(x_{3}\right)\right]
\]
\[
+x_{2}\sin\left(x_{3}\right)\left[m\cos\left(x_{3}\right)\left(mu_{2}-gm\cos\left(x_{3}\right)\right)-Mgm\right]+Imu_{3}{\cos^2\left(x_{3}\right)}\Bigg].
\]
\begin{rmk}\label{approx}
The use of the control force $F_{\theta},$ represented by control variable $u_3,$ admits the greater angles of the payload sway during transportation, and therefore, one must work with complete system, and the often used small-angle approximations of the type $\sin\theta\approx 0$ (or $\sin\theta\approx\theta$), $\dot\theta^2\approx 0$  and $\cos\theta\approx 1$ can not be applied in our analysis.
\end{rmk}
\section{Theoretical background. Control law design}

Our approach to the asymptotic stabilization of the overhead crane system is based on the two cornerstones of modern control theory and theory of dynamical systems.

Consider a linear time-invariant (LTI) system $\dot x=Ax+Bu,$ where $A$ and $B$ are $n\times n$ and $n\times m$ constant real matrices, respectively. A fundamental result of linear control theory is that the following three conditions are equivalent, see, e.~g. \cite{AntMich}: 
\begin{itemize}
\item[(i)] the pair $(A,B)$ is controllable;
\item[(ii)] $\rank\mathcal{C}_{(A,B)}=n,$ where $\mathcal{C}_{(A,B)}=:\left( B\  AB\  A^2B\ \cdots\ A^{n-1}B\right)$ is an $n\times mn$ Kalman's controllability matrix;
\item[(iii)] for every $n$-tuple real and/or complex conjugate numbers $\lambda_1,$ $\lambda_2,$ $\dots,$ $\lambda_n,$ there exists an $m\times n$ state feedback gain matrix $K$ such that the eigenvalues of the closed-loop system matrix $A_{cl}=A-BK$
are the desired values $\lambda_1,\lambda_2,\dots,\lambda_n.$
\end{itemize}

In general, the nonlinear control system
\begin{equation*}
\dot x=G(x,u), \ t\geq 0,\ x\in\mathbb{R}^n\ u\in\mathbb{R}^m, \ \dot{} =:d/dt,
\end{equation*}
with the state feedback of the form $u=-Kx$ and for which is assumed that $x=0$ is its solution, that is, $G(0,0)=0,$ is studied. It is well-known, that if the pair $(A,B),$ where $A=G_x(0,0)$  and $B=G_u(0,0)$ are the corresponding Jacobian matrices with respect to the state and input variables, respectively, and evaluated at $(0,0)$  is controllable, then the LTI system $\dot x=(A-BK)x$ is in some neighborhood of the origin topologically equivalent, and preserving the parametrization by time, to the system $\dot x=G(x,-Kx),$ provided that the eigenvalues of the matrix $A-BK$ have non-zero real part. The precise statement about this property gives the Hartman-Grobman theorem, see, e.~g. \cite[p.~120]{Perko},  providing the exact geometric characterization of the trajectories of the closed-loop system in the neighborhood of the equilibrium state.  Thus, if the matrix $K$ is chosen such that all eigenvalues of $A-BK$ have negative real parts, the nonlinear system $\dot x=G(x,-Kx)$ is locally asymptotically stable in the neighborhood of $x=0.$

Here is meant the usual definition of local asymptotic stability, that is, the solution $x=0$ of the system $\dot x=G(x,-Kx)$ is asymptotically stable if for every $\varepsilon>0$ there exists $\delta>0$ such that if $||x(0)||\leq\delta$ then $||x(t)||\leq\varepsilon$ for all $t\geq 0,$ and, moreover, $||x(t)||\rightarrow 0^+$ with $t\rightarrow \infty,$  see, e.~g., \cite[p.~19]{Barbashin}. Here $||\cdot||$ denotes the Euclidean vector norm.

Now, by imposing the natural requirement for the real working devices on the boundedness of the state variables $x_i,$ $i=1,\dots,6,$ namely that $||x||\leq\Delta$ for some constant $\Delta>0,$ the local asymptotic stability of the closed-loop system $\dot x=G(x,-Kx)$ around its equilibrium position $x=0$ by (formal and general) computing the lower bound of the region of attraction will be proved. In the engineering practice, one do not need to know this region explicitly because all generated trajectories are verified and validated for their suitability and appropriateness, and, moreover, a theoretical analysis of such highly coupled control systems that are investigated in the present paper is practically impossible. 

Let the gain matrix $K$ is such that the real parts of all eigenvalues of the matrix $A_{cl}$ are  negative. Then 
\[
\dot x=G(x,-Kx)=A_{cl}x+R_1(x),
\]
where $R_1(x)$ is the Taylor's remainder, obviously $||R_1(x)||=o(||x||)$ as $||x||\rightarrow 0^+.$ This implies that there exists a constant $\sigma=\sigma(K)>0$ such that $||R_1(x)||\leq\sigma(K)||x||$ for all $||x||\leq\tilde\Delta(\sigma).$  Let us consider as a Lyapunov function candidate $V(x)=x^TPx,$ where symmetric and  positive definite matrix $P$ is a solution of Lyapunov equation 
\[
PA_{cl}+A^T_{cl}P=-Q(K)
\]
for appropriate choice of the symmetric and positive definite matrix $Q,$ which can be solved as an optimization problem with regards to the gain matrix $K.$ Let the constant $\sigma(K)$ is such that 
\[
\lambda_{\min}(Q(K))>2\sigma(K)\lambda_{\max}(P),
\]
where $\lambda_{\min}$ and $\lambda_{\max}$ denote the minimal and maximal eigenvalue of the matrix, respectively.
Then along the trajectories of the system $\dot x=G(x,-Kx)=:g_K(x)$ is
\[
\dot V(x(t))=x^T(t)Pg_K(x(t))+g^T_K(x(t))Px(t)
\]
\[
=x^TP[A_{cl}x+R_1(x)]+[x^TA^T_{cl}+R^T_1(x)]Px
\]
\[
=x^T(PA_{cl}+A^T_{cl}P)x+2x^TPR_1(x)=-x^TQx+2x^TPR_1(x),
\]
by using the fact that $x^TPR_1(x)$ is a scalar, that is,  $x^TPR_1(x)=\left(x^TPR_1(x)\right)^T=R^T_1(x)Px.$ 

Because for each $n\times n$ symmetric and positive definite real matrix $C$ is
\[
\lambda_{\min}(C)||x||^2\leq x^TCx\leq\lambda_{\max}(C)||x||^2,\ x\in\mathbb{R}^n
\]
(a special case of Rayleigh-Ritz's theorem, \cite[p.~176]{HornJohnson}) and 
\[
x^TPR_1(x)=\frac{x^TPx}{||x||^2}x^TR_1(x)\leq\lambda_{\max}(P)||x||||R_1(x)||\leq\sigma(K)\lambda_{\max}(P)||x||^2,
\]
one get that
\[
\dot V(x(t))\leq-\bigg( \lambda_{\min}(Q(K))-2\sigma(K)\lambda_{\max}(P) \bigg)||x||^2<0
\]
for $||x||\leq\tilde\Delta,$ $x\neq0,$ which implies local asymptotic stability of the zero solution of the closed-loop system $\dot x=G(x,-Kx).$

One of the main aim of the present paper lies in comparing the performance of designed simulation- and state feedback-based feed-forward control for overhead crane system stabilizing the system in its end position for the crane with variable and constant length of the hoisting rope, Section~\ref{variable_length} and Section~\ref{constant_length}, respectively.
\section{Application to the overhead crane with variable length of rope. Simulation experiment in MATLAB}\label{variable_length}
First, let us verify that linear part of the overhead crane system $\dot x =G(x,u)$ is controllable at the point $(x,u)=(0,0):$  
\[
A=G_x(0,0)=
\left(\begin{array}{cccccc} 0 & 0 & 0 & 1 & 0 & 0\\ 0 & 0 & 0 & 0 & 1 & 0\\ 0 & 0 & 0 & 0 & 0 & 1\\ 0 & 0 & 0 & 0 & 0 & 0\\ 0 & 0 & 0 & 0 & 0 & 0\\ 0 & 0 & 0 & 0 & 0 & 0 \end{array}\right),\
B=G_u(0,0)=
\left(\begin{array}{ccc} 0 & 0 & 0\\ 0 & 0 & 0\\ 0 & 0 & 0\\ 1 & 0 & 0\\ 0 & 1 & 0\\ 0 & 0 & 1 \end{array}\right).
\]
Because $A^2=0,$ for the rank of the controllability matrix the equality 
\[
\rank \mathcal{C}_{A,B}=\rank (B\ AB\ A^2B\  A^3B\  A^4B\ A^5B)=\rank (B\ AB)
\]
holds and since
\[
(B\ AB)=
\left(\begin{array}{cccccc} 
0 & 0 & 0 & 1 & 0 & 0\\ 
0 & 0 & 0 & 0 & 1 & 0\\
 0 & 0 & 0 & 0 & 0 & 1\\ 
1 & 0 & 0 & 0 & 0 & 0\\ 
0 & 1 & 0 & 0 & 0 & 0\\ 
0 & 0 & 1 & 0 & 0 & 0 
\end{array}\right)
\]
with  $\det(B\ AB)= -1\neq0,$ the linear part of the system $\dot x=G(x,u)$ is controllable. As has been analyzed above, all eigenvalues of the closed-loop matrix $A_{cl}=A-BK$ can be arbitrarily assigned by appropriately selecting a state feedback gain matrix $K.$
\begin{rmk}
At this place it is worth noting that if the sway angle control force $F_{\theta}$ is not considered, the crane system may not be locally asymptotically stabilizable at the desired final position by using the linear state feedback control law.
\end{rmk}
Let the desired (and permissible) eigenvalues of the closed-loop system $\dot x=A_{cl}x$ are $p=[-0.1\ -0.15\ -0.2\ -0.25\ -0.3\ -0.35]$ that are determined in practice on the basis of the crane work parameters such as maximum permissible velocity of the used equipments, for example. In general, it is advisable to choose these eigenvalues real to avoid the oscillating trajectories. For the (negative) real eigenvalues, the convergence to the equilibrium point for $t\rightarrow \infty$ will be monotonous. The MATLAB output of the pole placement command {\tt K = place(A,B,p)} gives 
\[
K= 
\left(\begin{array}{ccccccc} 
0.1050     &    0    &     0  &  0.6500    &     0    &     0 \\
     0      &   0.0300   &      0   &      0  &  0.3500   &      0\\
         0    &     0   & 0.0250   &      0   &      0   & 0.3500

\end{array}\right),
\]
which locally asymptotically stabilizes the equilibrium point $x_e=0$  of the system $\dot x=G(x,-Kx).$
The routine {\tt place} in MATLAB uses the algorithm of \cite{KaNiDo} which, for multi-input systems, optimizes the choice of eigenvectors for a robust solution  and the sensitivity of the assigned poles to perturbations in the system and gain matrices is minimized. A general theory regarding pole placement problem for linear systems can be found in the work \cite[p.~335]{AntMich}.

Now, let the desired end position of the payload is $\tilde x_{e,1} = 10,$ $\tilde x_{e,2} = 3$ and $\tilde x_{e,i} = 0,$ for $i=3,4,5,6,$ corresponding to the placement $z=\tilde x_{e,1}$ and $y=-\tilde x_{e,2}$ in Fig.~\ref{fig:M1}. The state variables transformation $\tilde x=\Phi (x)$ defined by the formula
\begin{equation*}
\Phi: \tilde x_1= \tilde x_{e,1}+x_1,\, \tilde x_2=\tilde x_{e,2}-x_2,\, \tilde x_3=x_3,\, \tilde x_4 =x_4, \tilde x_5 =-x_5,\, \tilde x_6=x_6
\end{equation*}
is used. Obviously, $\Phi^{-1}(\tilde x_e)=0.$   Then the original system $\dot x=G(x,u)$ transformed to the new equilibrium point $\tilde x_e$ is $\dot{\tilde x}=\tilde G\left(\tilde x, \tilde u\right),$ where
\[
\tilde G\left(\tilde x, \tilde u\right)=\left(\tilde x_4,\tilde x_5,\tilde x_6,G_4\left(\Phi^{-1}(\tilde x), \tilde u\right),-G_5\left(\Phi^{-1}(\tilde x),\tilde u\right), G_6\left(\Phi^{-1}(\tilde x), \tilde u\right) \right)^T,
\]
$\tilde u$ is instead of $u.$ The direct calculation gives that
$\tilde A=:\tilde G_{\tilde x}(\tilde x_e,0)=A$ and 
\[
\tilde B=:\tilde G_{\tilde u}(\tilde x_e,0)=
\left(
\begin{array}{ccc} 
0 & 0 & 0\\ 0 & 0 & 0\\ 0 & 0 & 0\\ 1 & 0 & 0\\ 0 & -1 & 0\\ 0 & 0 & 1 
\end{array}
\right).
\]
It will be finding the feedback control law in the form $\tilde u=-\tilde Kx,$ and so  $\dot{\tilde x}=\tilde G\left(\tilde x, -\tilde K\Phi^{-1}(\tilde x)\right)$ with $\tilde A_{cl}=\tilde A-\tilde B\tilde K\Phi_{\tilde x}^{-1}(\tilde x_e),$ where $\Phi_{\tilde x}^{-1}(\tilde x_e)$ denotes the Jacobian matrix of the inverse of the transformation $\Phi$ and evaluated at $\tilde x_e.$ The sufficient condition to have the matrices $A_{cl}=A-BK$ and $\tilde A_{cl}$ the same eigenvalues is to be $BK=\tilde B\tilde K\Phi_{\tilde x}^{-1}(\tilde x_e),$ that is, the feedback gain matrix $\tilde K$  may be directly derived from the matrix $K$ by multiplying the second row and second and fifth column of $K$ by the number $(-1).$

For the purpose of numerical simulation the following data will be used: 
\begin{itemize}
\item[] $M=0.2\ [\times10^3\,\mathrm{kg}],\ $ $m=10\ [\times10^3\,\mathrm{kg}],\ $ $I=4\ [\times10^3\,\mathrm{kg\, m^2}],\ $ $g = 9.81\ [\mathrm{m\, s^{-2}}].$ 
\end{itemize}

\begin{figure} 
   \centerline{
    \hbox{
     \psfig{file=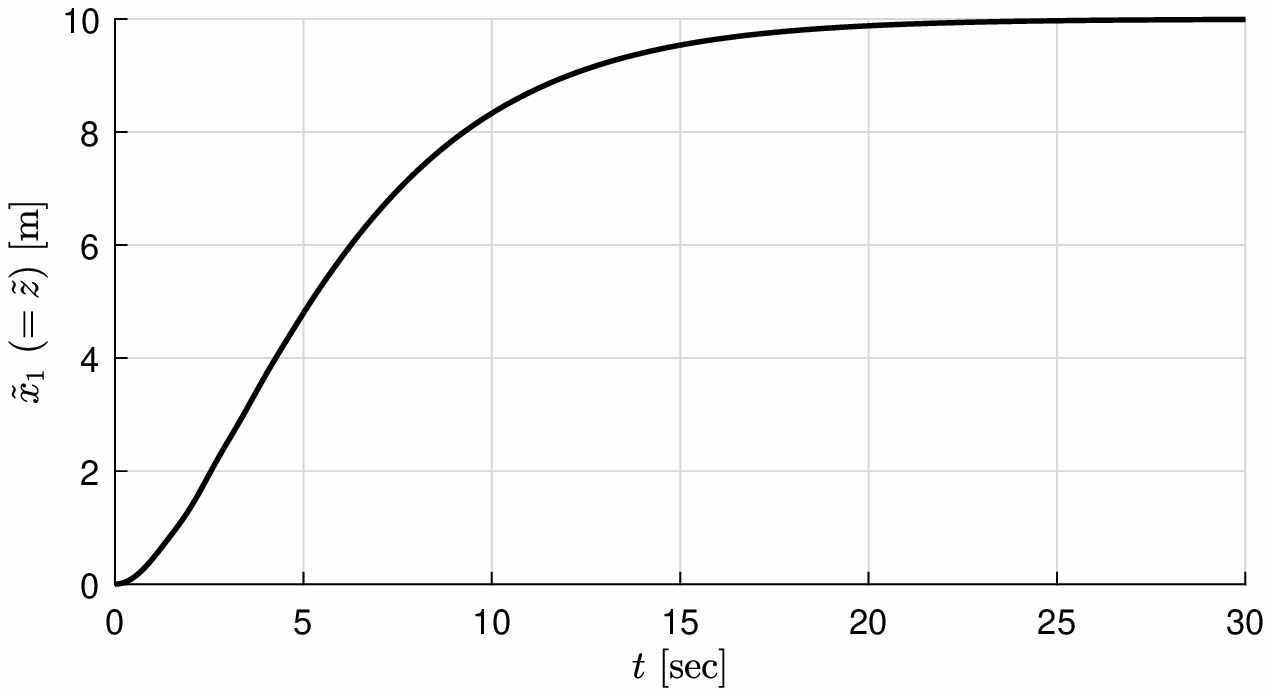,width=5.cm, clip=}
     \hspace{1.cm}
     \psfig{file=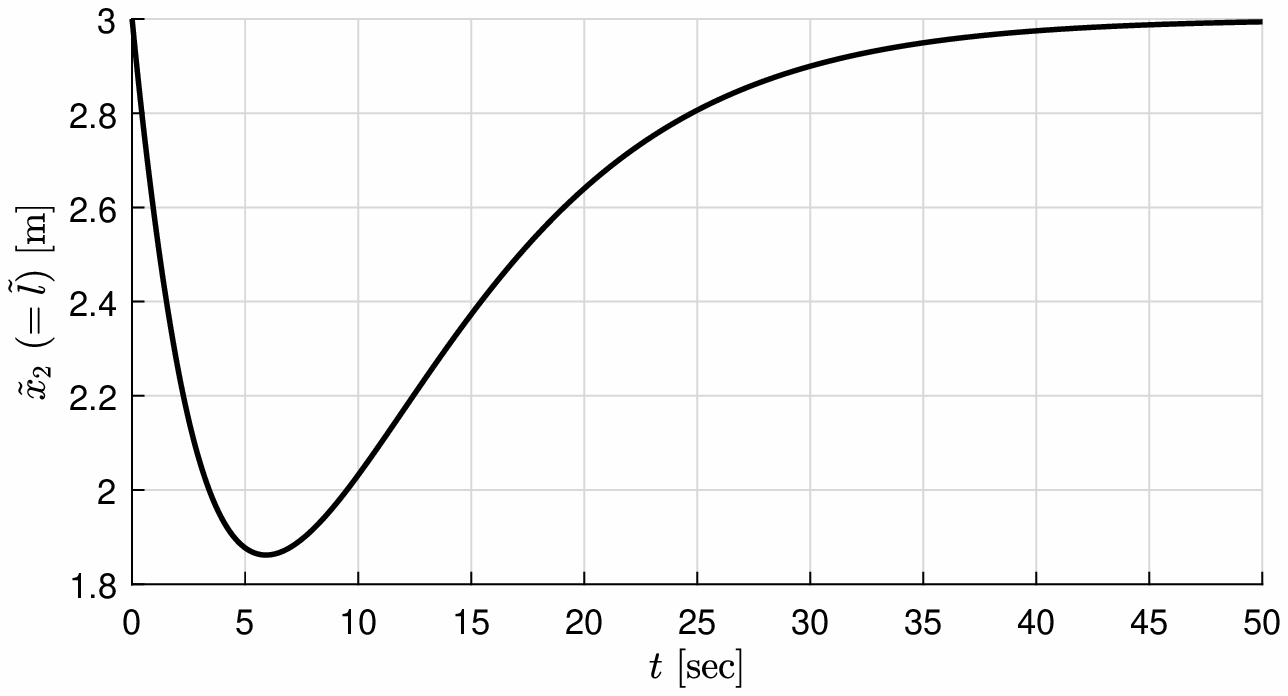,width=5cm,clip=}
    }
   }
   \vspace{0.5cm}
   \centerline{
    \hbox{
     \psfig{file=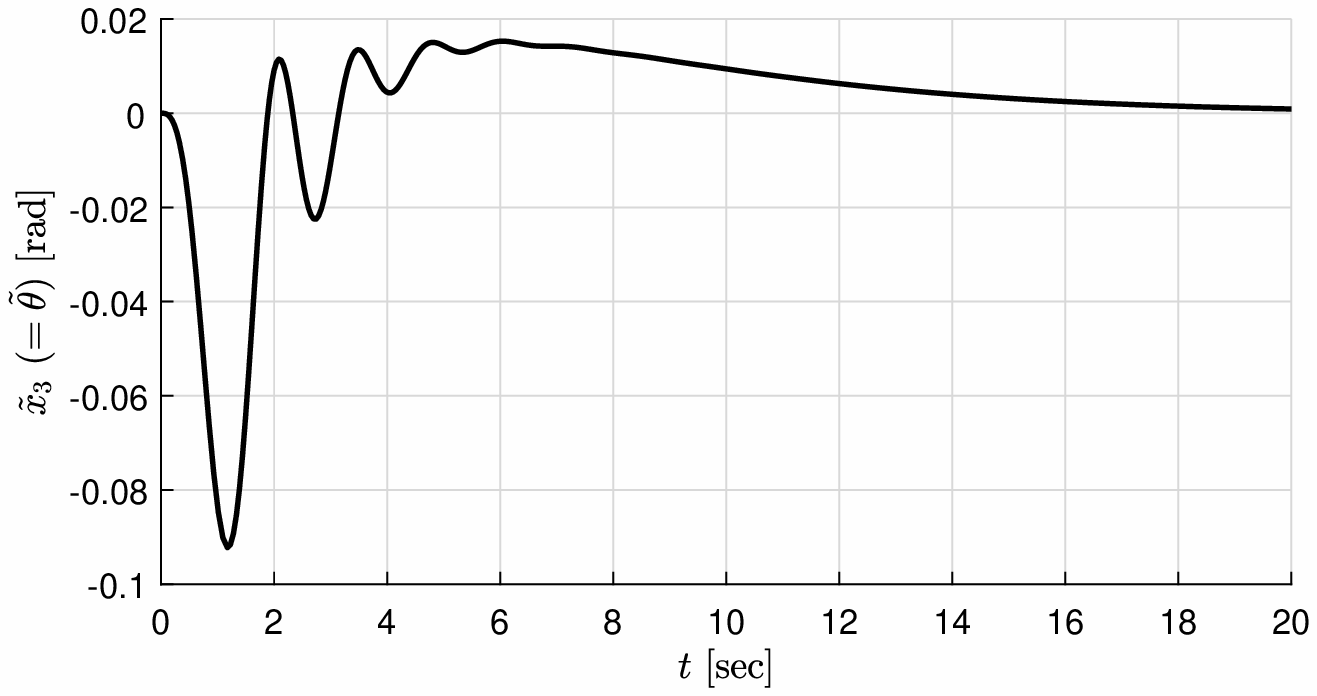,width=5.cm,clip=}
     \hspace{1.cm}
     \psfig{file=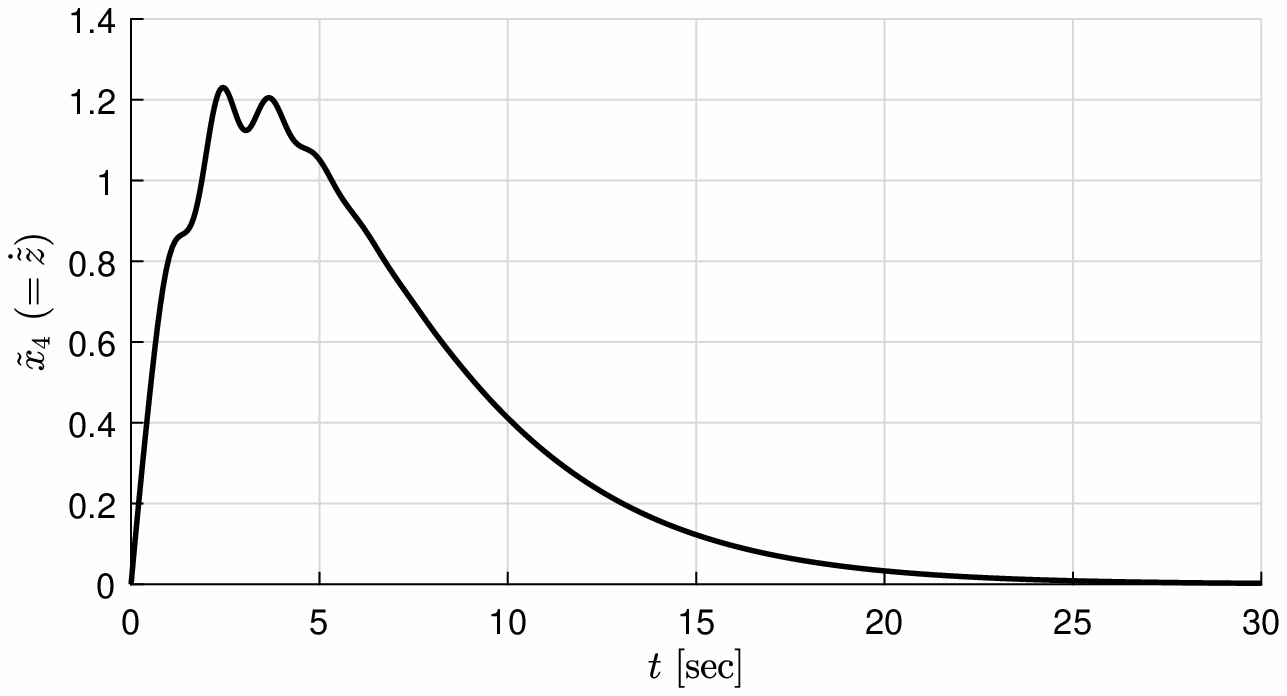,width=5.cm,clip=}
    }
   }
   \vspace{0.5cm}
   \centerline{
    \hbox{
     \psfig{file=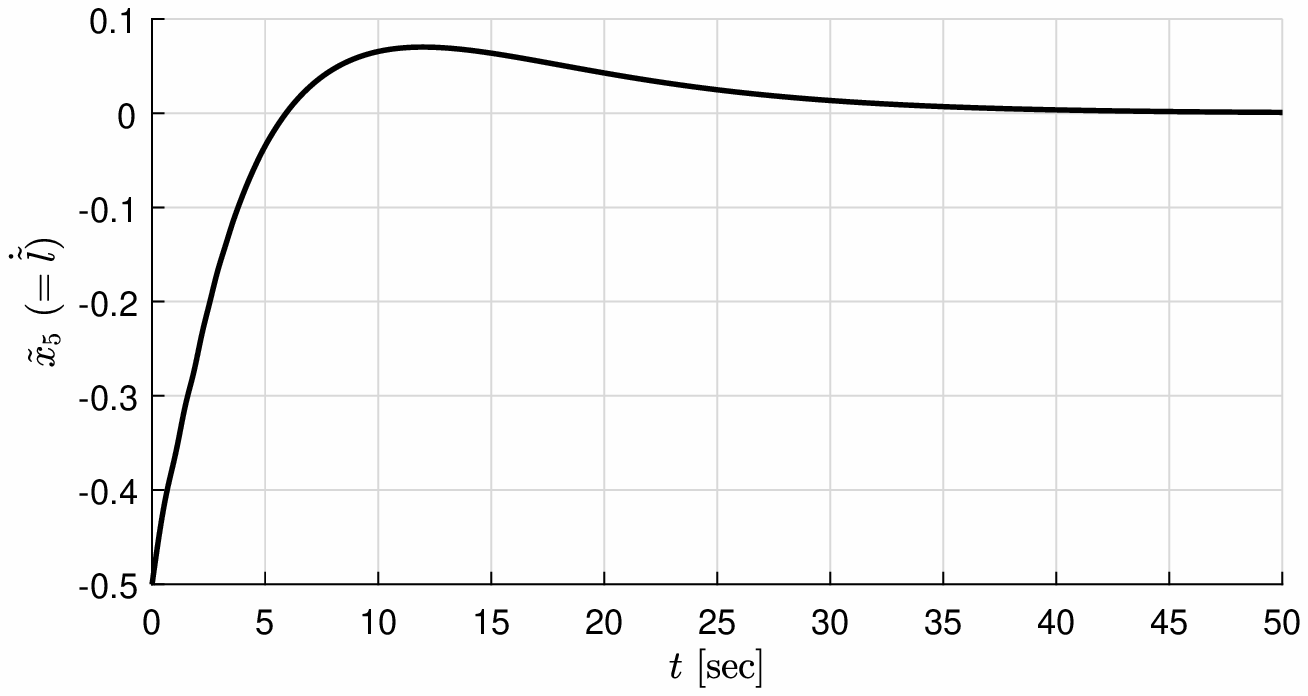,width=5.cm,clip=}
     \hspace{1.cm}
     \psfig{file=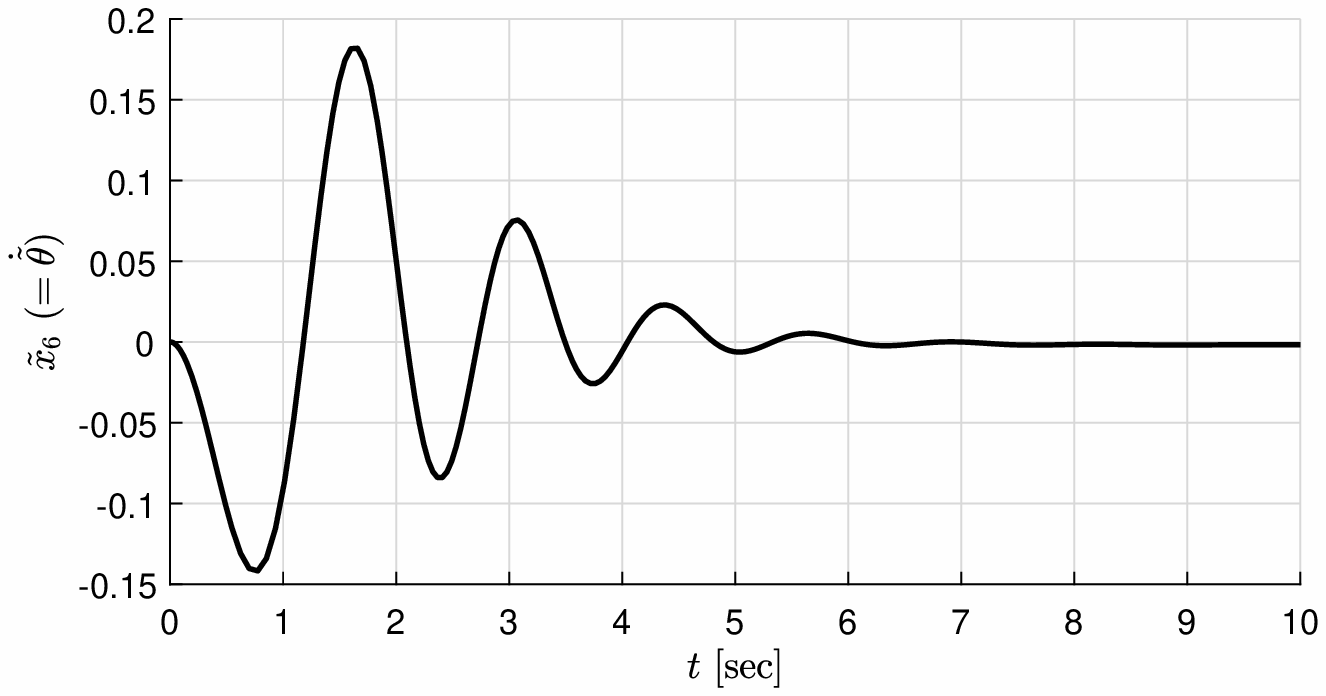,width=5.cm,clip=}
         }
   }
\caption{The solutions $\tilde x_i,$ $i=1,\dots,6$ of the reference model with varying rope length}
\label{solutions_xi}
\end{figure}

\begin{figure}
   \centerline{
    \hbox{
     \psfig{file=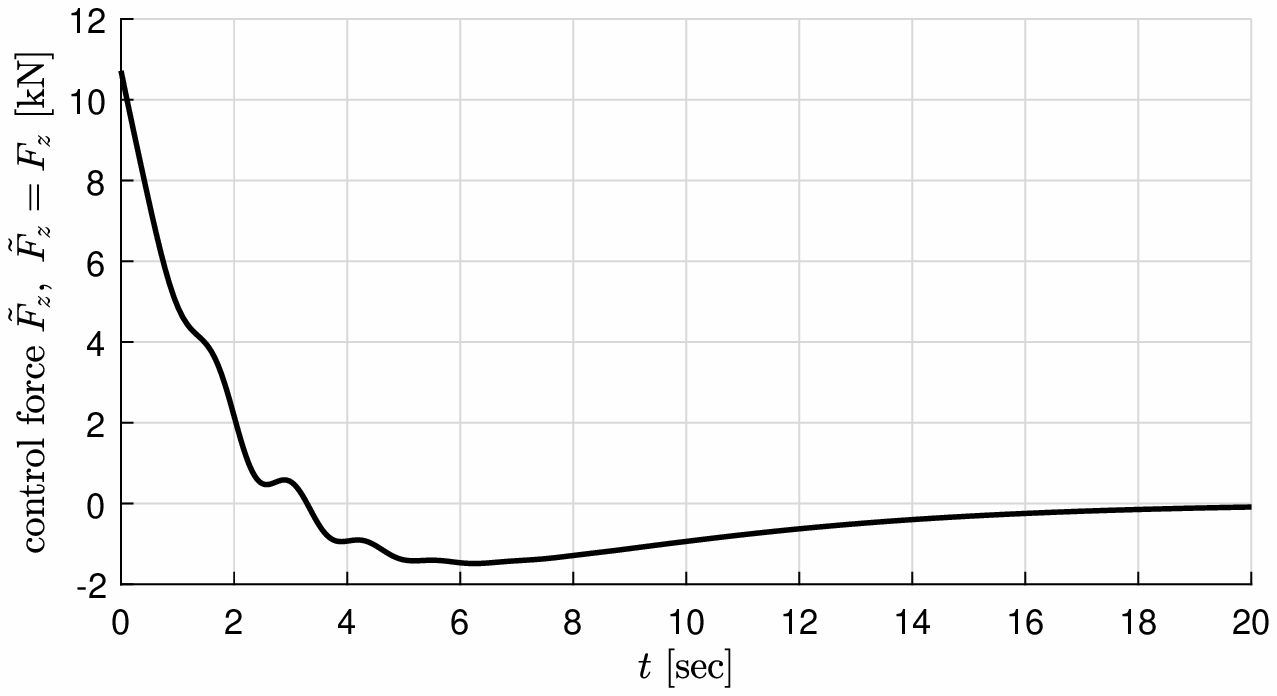,width=5cm, clip=}
     \hspace{1.cm}
     \psfig{file=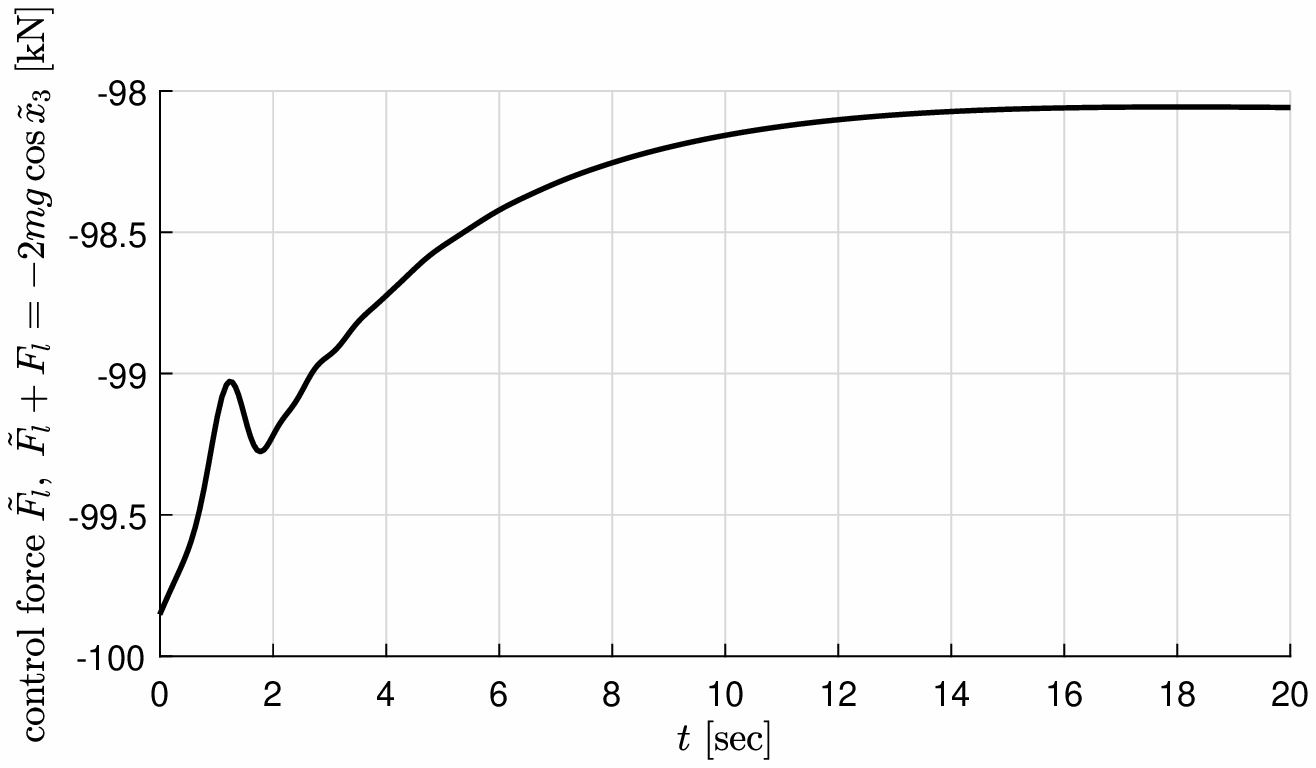,width=5cm,clip=}
    }
    \hbox{
     \psfig{file=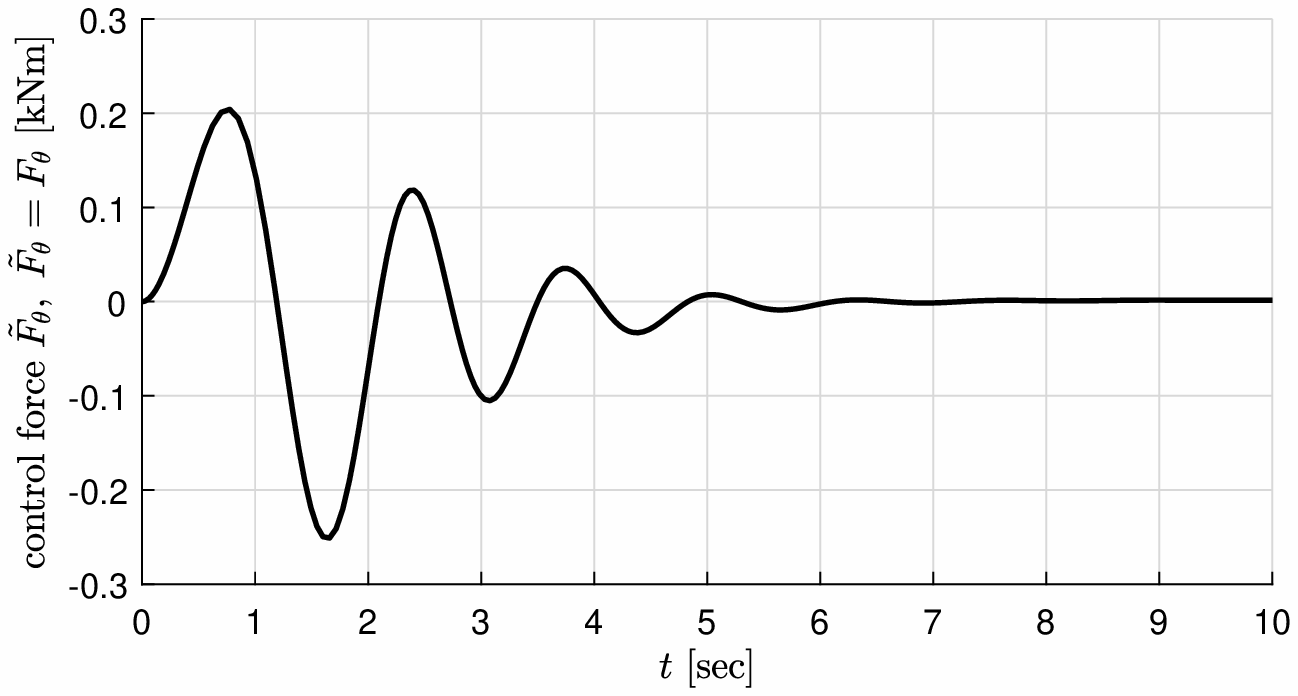,width=5cm, clip=} 
         }
   }   
\caption{The control forces $\tilde F_{z},$ $\tilde F_l$ and $\tilde F_{\theta}.$ The forces $\tilde F_{z}$ and $F_{z},$ and also $\tilde F_{\theta}$ and $F_{\theta}$ have equal magnitude and direction, $\tilde F_l$ and $F_l$ have opposite direction and $\tilde F_l+F_l=-2gm\cos(x_3)$}
\label{control_forces}
\end{figure}

Now, let the starting position is $ \tilde x(0) = (0\ 3\ 0\ 0\ -0.5\ 0)^T,$ corresponding $z=0$ and $y=-3$ in Fig.~\ref{fig:M1}, that is, the payload starts from rest and is pulled upwards with the velocity of $\tilde x_5=-0.5\ [\mathrm{m\, s^{-1}}].$ The time evolution of the state variables $\tilde x_i,$ $i=1,\dots,6$ is displayed in Fig.~\ref{solutions_xi} and the corresponding control forces $\tilde F_{z},$ $\tilde F_l$ and $\tilde F_{\theta}$  are depicted in Fig.~\ref{control_forces}.

Summarizing, the general strategy in stabilizing the system under consideration at the desired end position is as follows:
\begin{itemize}
\item[i)]  First, the controllability of the linear part of mathematical model of the overhead crane system derived above must be verified;  
\item[ii)] Secondly, the range of the eigenvalues $\lambda_i$ of the closed-loop system taking into account the technical limitations of the specific overhead crane is established, for example, the maximum permissible velocities of the trolley and hoist device or the jerk during starting;
\item[iii)]  Thirdly, using the appropriate state variables transformation the equilibrium point $x=0$ of the system $\dot x=G(x,-Kx)$ is translated to the new position $\tilde x_e,$ which is the desired end position of payload;
\item[iv)] Fourthly, for the desired eigenvalues of the closed-loop system, the gain matrix $K$ using the MATLAB command {\tt place} is calculated;
\item[v)]  Finally, the numerical simulation in the MATLAB environment is performed and the appropriateness of the generated trajectory is verified. Subsequently, the control forces 
\[
\tilde F_{z}(t)=(M+m)\tilde u_1(t)\ [=F_{z}(t)], 
\]
\[
\tilde F_l(t)=m\tilde u_2(t)-gm\cos\left(x_{3}(t)\right)\ [\tilde F_l(t)+F_l(t)=-2gm\cos(x_3(t))], 
\]
and
\[
\tilde F_{\theta}(t)=I\tilde u_3(t)\ [=F_{\theta}(t)],
\]
$\tilde u(t)=-\tilde Kx(t),$ stabilizing the overhead crane system at the desired end location of payload, to the  system are applied.
\end{itemize}

\newpage

\section{Application to the overhead crane with constant length of rope. Simulation experiment in MATLAB}\label{constant_length}
For comparison purpose, in this section the simulation experiment with the same data as in previous section is performed, with this difference that $\tilde l(t)=\tilde x_2(t)\equiv \tilde x_{e,2}=3$  and so $\dot{\tilde l}(t)=\tilde x_5(t)\equiv0.$ Thus, the number of state variables and governing equations reduces to four and only one control force is considered, $F_z,$ so, the control system is underactuated. Specifically, from the equations (\ref{eq:z}) and (\ref{eq:theta}) with $F_{\theta}\equiv0,$  one get for $\dot x_4$ and $\dot x_6$ the system of linear equation
\[
\left(
\begin{array}{cc} 
M+m & lm\cos\left(x_{3}\right)\\ lm\cos\left(x_{3}\right) & ml^2+I 
\end{array}
\right)
\left(
\begin{array}{c}
\dot x_4\\
\dot x_6
\end{array}
\right)
\]
\[
=\left(
\begin{array}{c} 
F_{z}+lm{x^2_{6}}\sin\left(x_{3}\right)\\ 
-glm\sin\left(x_{3}\right) 
\end{array}
\right),
\]
and after substituting  $u_1$ for $F_{z}$ one get the system  $\dot x =G(x,u),$ $x=(x_1,x_3,x_4,x_6)^T,$ $u=u_1,$ $G=(G_1,G_3,G_4,G_6)^T,$ where
\[
G_1=x_4,
\]
\[
G_3=x_6,
\] 
\[
G_4 =\frac{1}{l^2m^2{\sin^2\left(x_{3}\right)}+Ml^2m+Im+IM} 
\]
\[
\times\Bigg[Iu_{1}+\sin\left(x_{3}\right)\left(l^3m^2{x^2_{6}}+Ilm{x^2_{6}}\right)
\]
\[
+l^2mu_{1}+\frac12gl^2m^2\sin\left(2x_{3}\right)\Bigg],
\]
\[
G_6 = -\frac{1}{l^2m^2{\sin^2\left(x_{3}\right)}+Ml^2m+Im+IM}
\]
\[
\times\Bigg[\cos\left(x_{3}\right)\left[l^2m^2{x^2_{6}}\sin\left(x_{3}\right)+lmu_{1}\right]
\]
\[
+lm\left(Mg+gm\right)\sin\left(x_{3}\right)\Bigg].
\]
Now, the controllability of linear part of this system with
\[
A=G_x(0,0)= \left(\begin{array}{cccc} 0 & 0 & 1 & 0\\ 0 & 0 & 0 & 1\\ 0 & \frac{gl^2m^2}{Mml^2+I\left(M+m\right)} & 0 & 0\\ 0 & -\frac{glm\left(M+m\right)}{Mml^2+I\left(M+m\right)} & 0 & 0 \end{array}\right),
\]
\[
B=G_{u_1}(0,0)=\left(\begin{array}{c} 0\\ 0\\ \frac{ml^2+I}{Im+M\left(ml^2+I\right)}\\ -\frac{lm}{Mml^2+I\left(M+m\right)} \end{array}\right),
\]
will be verified. The controllability matrix $\mathcal{C}_{(A,B)}=(B\ AB\ A^2B\  A^3B)$ is a square matrix and 
\[
\det\mathcal{C}_{(A,B)} = -\frac{g^2 l^4 m^4}{{\left(M m l^2+I \left(M+m\right)\right)}^4}\neq 0,
\]
which implies the controllability of the linear part of the system under consideration.

\begin{figure}
   \centerline{
    \hbox{
     \psfig{file=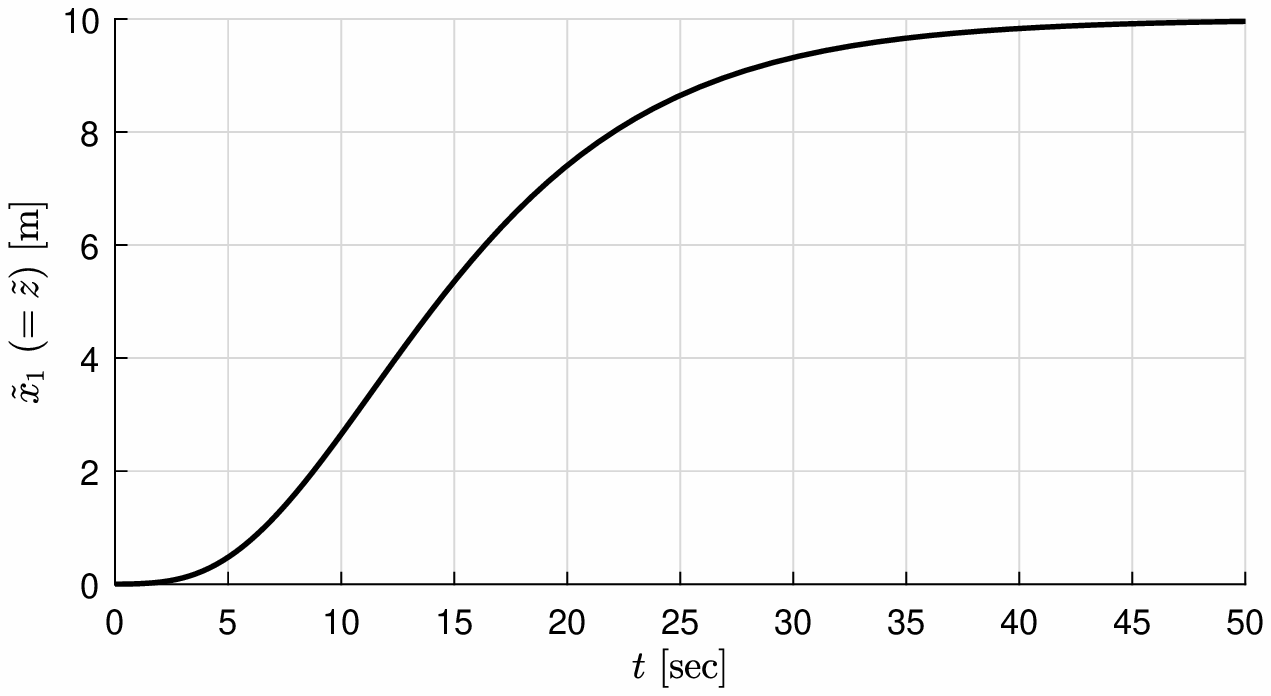,width=5.cm, clip=}
     \hspace{1.cm}
     \psfig{file=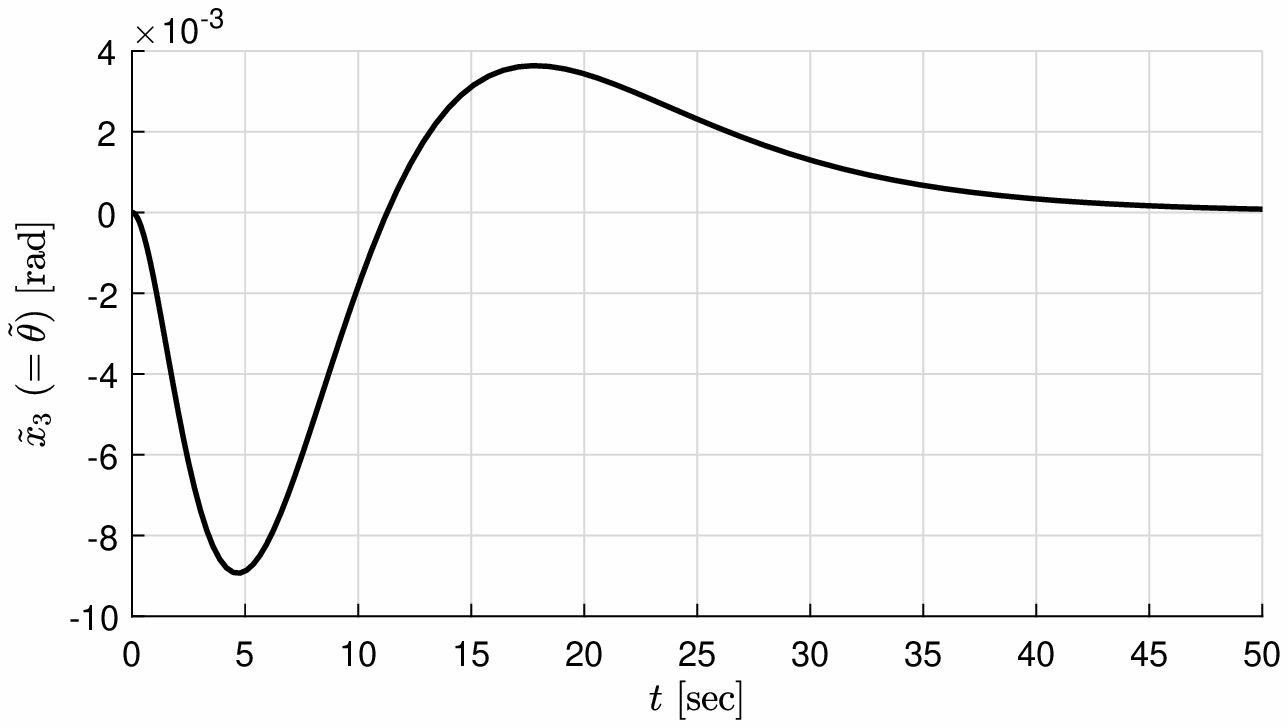,width=5cm,clip=}
    }
   }
     \centerline{
    \hbox{
     \psfig{file=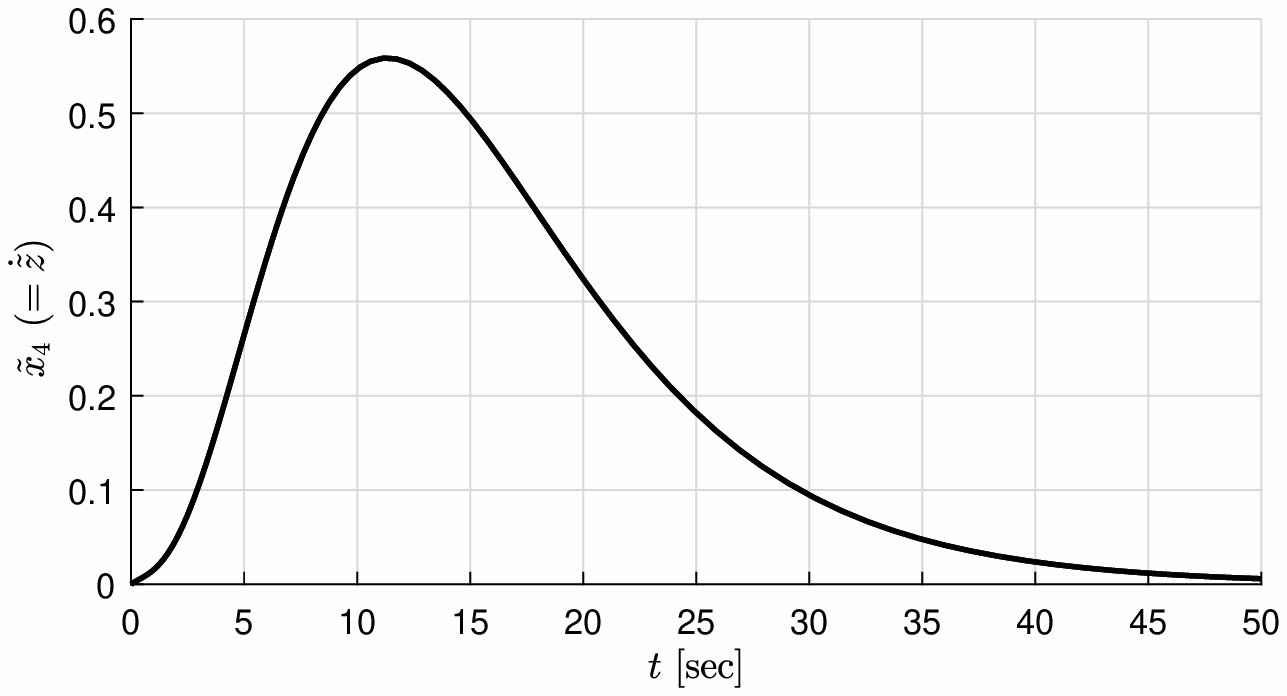,width=5.cm,clip=}
     \hspace{1.cm}
     \psfig{file=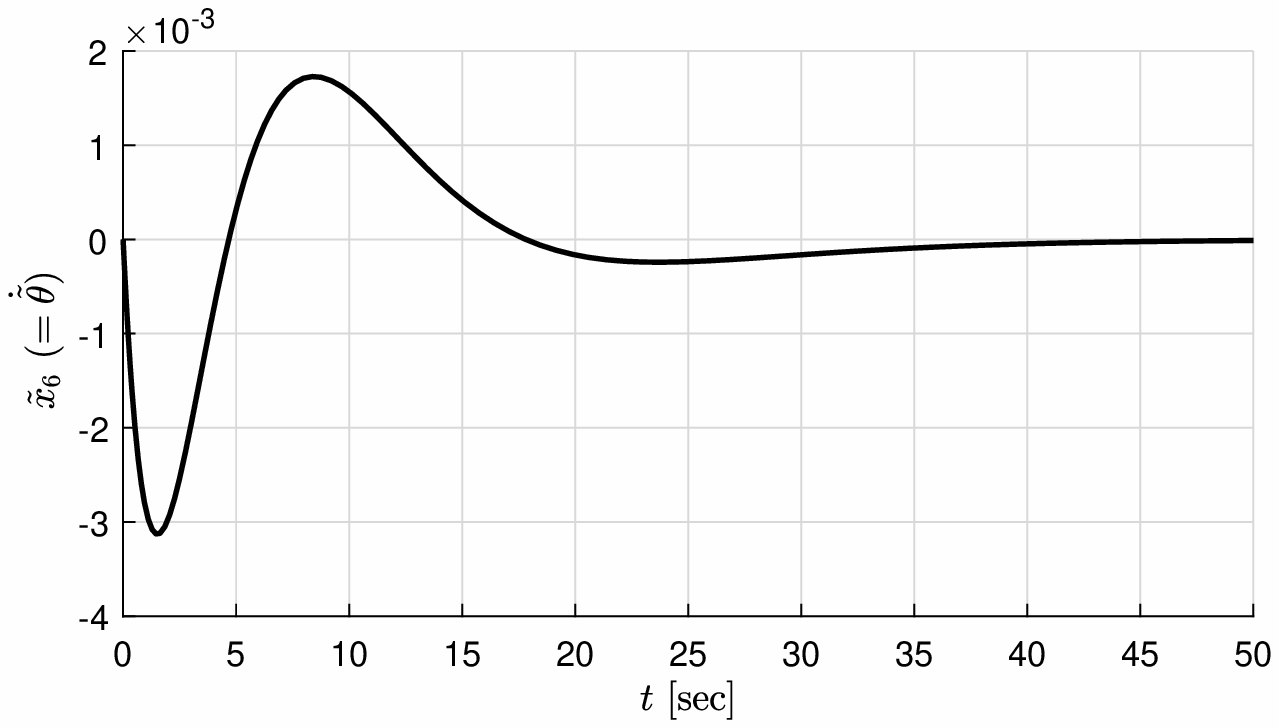,width=5.cm,clip=}
    }
   }     
\caption{The solutions $\tilde x_i,$ $i=1,3,4,6$ of the reference model with constant rope length}
\label{solutions_lconst_xi}
\end{figure} 
\begin{figure}[H]
   \centerline{
    \hbox{
     \psfig{file=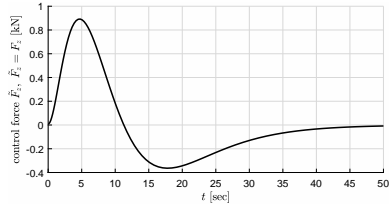,width=5cm, clip=}
         }
 }  
\caption{The control force $\tilde F_{z}$}
\label{control_force_Lconst}
\end{figure}  
For the simulation experiment the eigenvalues  
$p = [-0.2\ -0.25\ -0.3\ -0.35]$
of the closed-loop system are used, representing a conservative selection from the set of eigenvalues used in the previous section, and which are achieved for the gain matrix 
\[
K=(0.0010\   99.1882\    0.0159 \  -2.1061).
\]
The starting position is the same as in the previous simulation with a variable rope length, namely, 
$\tilde x(0)=(0\ 0\ 0\ 0)^T.$ In the Figs.~\ref{solutions_lconst_xi} and~\ref{control_force_Lconst} is depicted the time evolution of the state variables $\tilde x_1,$ $\tilde x_3,$ $\tilde x_4,$ $\tilde x_6$ and the control force $\tilde F_z,$ respectively.

Comparing the corresponding figures, one can see substantial prolongation of the transportation time to the end position although with significantly less payload sway.

\section{Conclusions}

In this paper, a linear control law for fully automated overhead crane systems with the aim to suppress the sway motion and to reduce the overall time of transportation using the state feedback-based feed-forward was proposed. It was shown that by the appropriate choice of the state feedback gain matrix the crane system can be asymptotically stabilized around the desired end position. Although the technical realization of the additional device for control of the sway angle (to ensure the controllability of linear part of system) requires some one-off costs of implementation, the numerical simulations indicate a substantial reduction of the transportation time (up to 50\%) in comparison with the overhead crane system with fixed rope length, as demonstrate the first sub-figures on top-left in the Figs.~\ref{solutions_xi} and \ref{solutions_lconst_xi}, and which may be desirable under certain circumstances.


\begin{thebibliography}{99}

\bibitem{Park}  
H.~Park, D.~Chwa, and K.-S.~Hong, A Feedback Linearization Control of Container Cranes: Varying Rope Length, International Journal of Control, Automation, and Systems, Vol.~5, No.~4, (2007), pp.~379--387. Permanent URL: 
\begin{verbatim}
http://www.ijcas.com/admin/paper/files/IJCAS_v5_n4_pp.379-387.pdf
\end{verbatim}

\bibitem{Tagawa_et_al} 
Y.~Tagawa, Y.~Mori, M.~Wada, E.~Kawajiri and K.~Nouzuka, Development of sensorless easy-to-use overhead crane system via simulation-based control, Journal of Physics: Conference Series 744 (2016).  DOI:10.1088/1742-6596/744/1/012017


\bibitem{Tagawa2_et_al} 
Y.~Tagawa, Y.~Hashizume and K.~Shimono, Simulation-based Control and its Application to a Crane System, ASME Proceedings, Control, Monitoring, and Energy Harvesting of Vibratory Systems: Active Vibration Control, 10 pages (2013). DOI:10.1115/DSCC2013-4031 

\bibitem{Yang}
J.H.~Yang, K.S.~Yang, Adaptive coupling control for overhead crane systems, Mechatronics,
Vol.~17, Iss.~2–3 (2007), pp.~143--152. DOI: 10.1016/j.mechatronics.2006.08.004
 
\bibitem{Tuan}
L.A.~Tuan, Design of Sliding Mode Controller for the 2D Motion of an Overhead Crane with Varying Cable Length, Journal of Automation and Control Engineering, Vol. 4, No. 3 (2016), pp.~181--188. DOI: 10.18178/joace.4.3.181-188

\bibitem{Tomczyk}
J.~Tomczyk, J.~Cink, and A.~Kosucki, Dynamics of an overhead crane under a wind disturbance condition, Automation in Construction, Vol.~42 (2014), pp.~100--111. DOI: 10.1016/j.autcon.2014.02.013

\bibitem{Sorensen}
K.L. Sorensen, W. Singhose, and S. Dickerson, A controller enabling precise positioning and sway reduction in bridge and gantry cranes, Control Engineering Practice, 15 (2007), pp.~825–-837.
DOI:10.1016/j.conengprac.2006.03.005

\bibitem{Zhang}
M. Zhang, X. Ma, X. Rong, X. Tian, and Y. Li, Error tracking control for underactuated overhead cranes against arbitrary initial payload swing angles, Mechanical Systems and Signal Processing, 84 (2017) pp.~268-–285.
DOI: 10.1016/j.ymssp.2016.07.028
  
\bibitem{AntMich}
P.J. Antsaklis, A.N. Michel, Linear Systems, 2006 Birkhauser Boston, 2nd Corrected Printing, Originally published by McGraw-Hill, Englewood Cliffs, NJ, 1997. DOI: 10.1007/0-8176-4435-0, Hardcover ISBN 978-0-8176-4434-5

\bibitem{Perko} 
L.~Perko, Differential Equations and Dynamical Systems (3rd Ed.). Texts in Applied Mathematics 7, New York:  Springer-Verlag, 2001. DOI: 10.1007/978-1-4613-0003-8, Hardcover ISBN 978-0-387-95116-4


\bibitem{Barbashin}
 E.A. Barbashin,  Introduction to the theory of stability,  Wolters-Noordhoff, Groningen, 1970. ISBN: 9001053459
 
\bibitem{HornJohnson} 
R.A. Horn, C.R. Johnson, Matrix analysis, Cambridge University Press, 1990. https://doi.org/10.1002/zamm.19870670330

\bibitem{KaNiDo}
J.~Kautsky, N.K.~Nichols, and P.~Van Dooren, Robust Pole Assignment in Linear State Feedback,  International Journal of Control, 41 (1985), pp. 1129--1155. https://doi.org/10.1080/0020718508961188 

\end{thebibliography}
\end{document}